\documentclass[12pt]{amsart}
\usepackage{amssymb}
\usepackage{amsthm}
\usepackage{amscd}

\usepackage{color}
\usepackage{hyperref}

\usepackage{tikz-cd}
\usepackage[utf8]{inputenc}
\usepackage[T1]{fontenc}   

\usepackage{mathpazo}
\usepackage[scaled=.95]{helvet}
\usepackage{courier}

\swapnumbers
\headsep=1cm \numberwithin{equation}{section}
\usepackage[colorinlistoftodos]{todonotes}

\usepackage[all]{xy}

\newtheorem{theorem}{Theorem}[section]

\newtheorem{lemma}[theorem]{Lemma}

\newtheorem{proposition}[theorem]{Proposition}

\newtheorem{corollary}[theorem]{Corollary}

\theoremstyle{definition}

\newtheorem{definition}[theorem]{Definition}

\newtheorem{example}[theorem]{Example}

\newtheorem{remark}[theorem]{Remark}

\newtheorem{question}[theorem]{Question}

\numberwithin{equation}{section}

\newtheorem{remark/questions}[theorem]{Remark and Questions}
\newtheorem{fact}[theorem]{Fact}
\newtheorem{facts}[theorem]{Facts}

\newtheorem*{agra}{Acknowledgment}

\newtheorem*{fund}{Funding}

 \long\def\alert#1{\smallskip\line{\hskip\parindent\vrule
\vbox{\advance\hsize-2\parindent\hrule\smallskip\parindent.4\parindent
  \narrower\noindent#1\smallskip\hrule}\vrule\hfill}\smallskip}

\newtheorem{newclaim}[theorem]{}

\DeclareMathOperator{\pd}{pd}

\DeclareMathOperator{\Hom}{Hom}
\DeclareMathOperator{\Ext}{Ext}

\DeclareMathOperator{\qpd}{qpd}

\DeclareMathOperator{\depth}{depth}

\newcommand\cdim{\operatorname{CI-dim}}

\def\dim{\mathop{\rm dim}}

\numberwithin{equation}{section}

\begin{document}

\dedicatory{}

\title[]{Ischebeck's formula, Grade and quasi-homological dimensions}

\author{V. H. Jorge-P\'erez}
\author{Paulo Martins}
\author{V. D. Mendoza-Rubio}
\address{Universidade de S{\~a}o Paulo -
ICMC, Caixa Postal 668, 13560-970, S{\~a}o Carlos-SP, Brazil}
\email{vhjperez@icmc.usp.br}

\address{Universidade de S{\~a}o Paulo -
ICMC, Caixa Postal 668, 13560-970, S{\~a}o Carlos-SP, Brazil}
\email{paulomartinsmtm@gmail.com}
\address{Universidade de S{\~a}o Paulo -
ICMC, Caixa Postal 668, 13560-970, S{\~a}o Carlos-SP, Brazil}
\email{vicdamenru@usp.br}

\keywords{quasi-projective dimension, quasi-injective dimension, grade, vanishing of Ext, Ischebeck's formula, quasi-perfect module, Gorenstein dimension}
\subjclass[2020]{13D05, 13D07,13D02}

\thanks{Corresponding author: Victor D. Mendoza Rubio}


\begin{abstract}

The quasi-projective dimension and quasi-injective dimension are recently introduced homological invariants that generalize the classical notions of projective dimension and injective dimension, respectively. For a local ring $R$ and finitely generated $R$-modules $M$ and $N$, we provide conditions involving quasi-homological dimensions where the equality 
$\sup \{i\geq 0: \Ext_R^i(M,N)\not=0\}=\depth R-\depth M$, which we call Ischebeck's formula, holds. One of the results in this direction generalizes a well-known result of Ischebeck concerning modules of finite injective dimension, considering the quasi-injective dimension. On the other hand, we establish an inequality relating the quasi-projective dimension of a finitely generated module to its grade and introduce the concept of a quasi-perfect module as a natural generalization of a perfect module. We prove some results for this new concept similar to the classical results. Additionally, we provide a formula for the grade of finitely generated modules with finite quasi-injective dimension over a local ring, as well as grade inequalities for modules of finite quasi-projective dimension. In our study, Cohen-Macaulayness criteria are also obtained.

\end{abstract}

\maketitle

\section{Introduction}
Throughout this paper, we assume that \( R \) is a commutative Noetherian ring and that all \( R \)-modules are finitely generated. The study of homological dimensions of modules is an important topic in commutative algebra, closely connected to the resolution of central problems in this area, 
as evidenced by several works such as \cite{ArayaYoshino,ArayaTakahashi2022,GhoshTakahashi2024,CelikbasIimaSadeghiTakahashi2018,GJT,Gheibi}. For any non-zero $R$-modules  $M$ and $N$, we define:
\begin{align*}
\operatorname{P}_R(M,N)&:=\sup \lbrace i \geq 0 : \operatorname{Ext}^i_R(M,N) \neq 0 \rbrace, \\ \operatorname{grade}(M,N)&:=\inf \lbrace i\geq 0: \Ext_R^i(M,N)\not=0 \rbrace.
\end{align*}
These values have been widely studied, especially when at least one of the modules $M $ or $N $ has some finite homological dimension. For instance, regarding the $\operatorname{P}_R(M,N)$, when $R$ is local the equality
\begin{equation}\label{Fisch}
    \operatorname{P}_R(M,N)=\depth R-\depth M
\end{equation}
 was proved in several cases involving some homological dimension.  One of the more known cases is when $N$ has finite injective dimension, due to Ischebeck \cite[2.6]{isch}, letting thus the authors motivated to refer to \eqref{Fisch} as the \textit{Ischebeck formula}. In the same paper, Ischebeck showed that \eqref{Fisch} is true when $\operatorname{pd}_R M < \infty$, where $\operatorname{pd}_R M$ denotes the projective dimension of the $R$-module $M$. Other cases where the equality has been proved include: (1) $M$ has finite complete intersection dimension and $\operatorname{P}_R(M,N)<\infty$, by Araya and Yoshino \cite[Theorem 4.2]{ArayaYoshino}. (2) $M$ has finite injective dimension and $N$ has finite Gorenstein injective dimension, due to Sazeedeh \cite[Theorem 2.10]{saz}. (3) $M$ has finite AB-dimension by Araya \cite[Lemma 2.5]{Araya2}.    (4) $M$ has finite redu\-cing complete intersection dimension with $\operatorname{depth} M \leq \operatorname{depth} R$ and $\operatorname{P}_R(M,N)<\infty$, by Celikbas et al.  \cite[Proposition 4.9]{Celikbaset}.

 On the other hand, the study of the grade involving certain homological dimensions can be evidenced in works like \cite{ArayaYoshino,Yoshida,Yassemi,GLCMHGD,new,Foxby1,khatami}. The study of \textit{perfect modules}, that is, modules satisfying $\operatorname{grade} M := \operatorname{grade}(M, R) = \operatorname{pd}_R M$, has received particular attention in this context. 
 This concept was generalized by Foxby \cite{Foxby1}, who replaced projective dimension with Gorenstein dimension, introducing a new class of modules called $G$--perfect modules. This new class has been further explored in works such as \cite{Yassemi, GhoshPuth}, among others.
 


Recently,  Gheibi, Jorgensen and Takahashi \cite{GJT} introduced a new homological dimension called quasi-projective dimension (qpd) which is a finer invariant than projective dimension, in the sense that for any $R$-module $M$, there is always an inequality $\operatorname{qpd}_R M \leq \operatorname{pd}_R M$ and the equality holds when $\operatorname{pd}_R M$ is finite. Following the work \cite{GJT}, Gheibi \cite{Gheibi} introduced another new homological dimension called quasi-injective dimension (qid) as a generalization of the injective dimension. We refer to these new dimensions as \textit{quasi-homological dimensions}.  The aim of this paper is to study the values $\operatorname{P}_R(M,N)$ and $\operatorname{grade}(M,N)$ when at least one between $M$ and $N$ has a finite quasi-homological dimension.

The following theorem presents our main results concerning the study of \( \operatorname{P}_R(-,-) \) and quasi-homological dimensions, and is motivated by {\cite[Theorems 6.2(4) and 6.12]{GJT} and \cite[Corollary 3.5]{Gheibi}}. In addition, it provides new cases in that the Ischebeck formula holds under certain conditions of the finiteness of quasi-homological dimension on  $M$ or $N$, obtaining thus not only improvements of \cite[Corollary 3.5]{Gheibi} and \cite[Lemma 3.5(1)]{MCMtensorproductsandvanishingofExtmodules}, but also a generalization of the well-known result of Ischebeck concerning to finitely generated modules of finite injective dimension \cite[2.6]{isch} mentioned previously.
\begin{theorem}[See Theorems \ref{teo:qproj}, \ref{g-dim} and \ref{isch}]\label{qw88}
    Let $R$ be a Noetherian local ring, and let $M$ and $N$ be non-zero $R$-modules. Suppose that $\operatorname{P}_R(M,N)<\infty$. Then the equality $$\operatorname{P}_R(M,N)=\operatorname{depth} R-\depth M$$
    holds in each one of the following cases:
    \begin{enumerate}
        \item $M$ has finite quasi-projective dimension. 
        \item  $M$ has finite Gorenstein dimension and $N$ has finite quasi-projective dimension.
        \item $N$ has finite quasi-injective dimension.  
    \end{enumerate}
\end{theorem}



On the other hand, in our study of the grade of two modules $M$ and $N$ and quasi-homological dimension, we give special attention to the case where $N=R$.  We prove a relation of inequality between the grade and quasi-projective dimension, similarly as occurs with the projective dimension and Gorenstein dimension. 

\begin{theorem}[See Theorem \ref{theo4.5}]\label{to12}
Let $M$ be a non-zero $R$-module. One then has
the inequality  $\operatorname{grade} M \leq \operatorname{qpd}_R M$.
\end{theorem}

It is worth noting that, unlike the case of projective dimension and Gorenstein dimension, the inequality in the theorem above does not follow trivially. The motivation for establishing Theorem \ref{to12} is to introduce the notion of quasi-perfect modules as a natural generalization of perfect modules.  More specifically, we define an $R$-module $M$ to be \textit{quasi-perfect} if $\operatorname{qpd}_R M<\infty$ and $\operatorname{grade} M=\operatorname{qpd}_R M$. This new definition allows us to generalize several well-known results for perfect modules, where Theorem \ref{to12} plays a crucial role. For instance, we generalize a well-known result stating that: A (finitely generated) $R$-module $M$ over a regular local ring is Cohen-Macaulay if and only if it is perfect (see e.g. \cite[Corollary 2.2.10]{BH}). Namely, we prove that an $R$-module with finite quasi-projective dimension is Cohen-Macaulay if it is quasi-perfect, and that the converse holds when $R$ is Cohen-Macaulay (see Proposition \ref{prop:eqv}). A ce\-lebrated application of (Peskine-Szpiro) intersection theorem is that: A local ring is Cohen-Macaulay if it admits a non-zero Cohen-Macaulay module of finite projective dimension (see \cite[Corollary 9.6.2]{BH}). However, the same does not work with quasi-projective dimension in place of projective dimension. As an application of our results, we obtain a Cohen-Macaulayness criterion for a local ring in terms of the existence of modules with finite quasi-projective dimension (see Corollary \ref{pros1}).

Concerning the study of grade and quasi-injective dimension, we obtain a formula for the grade of finitely generated modules with finite quasi-injective dimension over a local ring, which improves \cite[Theorem 3.7]{Gheibi}. In fact, this theorem says that when $R$ is local and  $M$ is a non-zero $R$-module of finite quasi-injective dimension, one then has the inequality $\dim M\leq \operatorname{depth} R$. However, this result does not describe the difference $\operatorname{depth} R - \dim M$. The formula we establish shows that this difference is $\operatorname{grade} M$.

\begin{theorem}[See Theorem \ref{teo:forinj}]\label{12aqs}
    Let $R$ be a local ring and let $M$ be a non-zero $R$-module with $\operatorname{qid}_R M <\infty$. Then $$\operatorname{dim} M=\depth R-\operatorname{grade} M.$$
\end{theorem}

The above theorem can be compared to \cite[Theorem 4.2]{new}. Motivated by Bass's Conjecture (see e.g. \cite[Corollary 9.6.2]{BH}), Takahashi \cite[Theorem 3.5(1)]{Takahashi} proved that a local ring \( R \) with a dualizing complex is Cohen-Macaulay if it admits a non-zero \( R \)-module of finite Gorenstein injective dimension and maximal Krull dimension. Moreover, Gheibi \cite[Corollary 3.8]{Gheibi} proved an analogous to Takahashi's result with quasi-injective dimension. As an application of Theorem \ref{12aqs}, we derive the following Cohen-Macaulayness criterion.

\begin{corollary}[See Corollary \ref{corol:criteria}]
Let $R$ be a local ring. If there exists a non-zero $R$-module $M$ such that $\operatorname{qid}_R M < \infty$ and $\dim R = \dim M + \operatorname{grade} M$, then $R$ is Cohen-Macaulay.
\end{corollary}

Now, we briefly describe the structure of this paper. In Section \ref{section2}, we provide de\-finitions, notations, and some results that are considered in this paper. In Section \ref{section3}, we prove items (1) and (2) of Theorem \ref{qw88} and explore some of their consequences. In Section \ref{section4}, we prove item (3) of Theorem \ref{qw88} obtaining a generalization of a celebrated result by Ischebeck. In Section \ref{section5}, we prove results concerning the grade and the quasi-homological dimensions, we prove Theorems \ref{to12} and \ref{12aqs}, introduce the definition of quasi-perfect module and explore some applications obtaining Cohen-Macaulayness criteria for rings and modules. Finally, in the last section, we apply some results of Section 5 to obtain some grade inequalities for modules with finite quasi-projective dimension, with Theorem \ref{genary} being the main result of this section.



\section{Conventions and Background}\label{section2}

In this section, we introduce fundamental definitions and facts, such as  quasi-projective dimension, quasi-injective dimension and Gorenstein dimension. These concepts will be essential throughout the rest of the paper. Unless otherwise specified, all modules will be considered over a commutative Noetherian ring  $R$, not necessarily local.

\begin{newclaim}
Let $R$ be a local ring and let $M$ be an $R$-module. Consider a minimal free resolution
\begin{align*}
\cdots \rightarrow F_i \xrightarrow{\varphi_i} F_{i-1} \rightarrow \cdots \rightarrow F_1 \xrightarrow{\varphi_1} F_0 \xrightarrow{\varphi_0} M \rightarrow 0
\end{align*}
of $M$. For $i \geq 1$, the $i$-\textit{syzygy} of $M$, denoted by $\Omega^i(M)$, is defined as the kernel of the map $\varphi_{i-1}$. When $i=0$, we set $\Omega^0(M)=M$. For $i \geq 0$, the modules $\Omega^i(M)$ are defined uniquely up to isomorphism.
\end{newclaim}

\begin{newclaim}\label{convetion}
In this work, we adopt the conventions that $\sup \emptyset = -\infty$ and $\inf \emptyset = \infty$. We also adopt the convention that the Krull dimension of the zero module is $-\infty$.

\end{newclaim}

\begin{newclaim} Let $M$ and $N$ be non-zero $R$-modules. Recall the following notations:
\begin{align*}
    \operatorname{P}_R(M,N)&=\sup \{ i \geq 0: \Ext_R^i(M,N)\not=0\},\\
    \operatorname{grade}(M,N)&=\inf\{i \geq 0: \Ext_R^i(M,N)\not=0\}.
\end{align*}
 Note that $\operatorname{grade}(M,N)$ can be infinite (e.g., let $\mathfrak{m} \neq \mathfrak{n}$ be two maximal ideals of $R$, $M=R/\mathfrak{m}$ and $N= R/\mathfrak{n}$).
\end{newclaim}

Let $M$ be an $R$-module and let $I$ be an ideal of $R$ such that $IM \neq M$. We denote by $\operatorname{depth}(I,M)$ the common length of a maximal $M$-regular sequence in $I$ (see \cite[Definition 1.2.6]{BH}). 


\begin{newclaim}\label{claim2.3}
Let $R$ be a local ring, and let $M$ and $N$ be non-zero $R$-modules. One  then has $\operatorname{grade}(M,N)<\infty,$ and the following chain of inequalities holds: 
\begin{align*}
0 \leq \operatorname{grade}(M,N) \leq \operatorname{P}_R(M,N) \leq \infty. 
\end{align*}
\end{newclaim}
\begin{proof}
We only need to show that $\operatorname{grade}(M,N)<\infty$. 
Set $I= \operatorname{ann}(M)$. Note that $I$ is a proper ideal of $R$, as $M \neq 0$. Since $N \neq 0 $, by Nakayama's Lemma, we have $IN \neq N$. Thus, by \cite[Proposition 1.2.10 (e)]{BH}, $\operatorname{grade}(M,N) = \inf \lbrace i \geq 0 : \operatorname{Ext}_R ^i(M,N) \neq 0 \rbrace = \operatorname{depth}(I,M)$, whence $\operatorname{grade}(M,N) < \infty$. 
\end{proof}




 

\begin{facts}\label{facts3} Let $M$ and $N$ be non-zero $R$-modules.
\begin{enumerate}
\item \cite[Theorem 2.1]{Yassemi} The following inequalities hold:
\begin{itemize}
    \item[(a)] \(\operatorname{depth} N - \dim M \leq \operatorname{grade}(M, N)\);
    \item[(b)] If \(\operatorname{Supp} M \subseteq \operatorname{Supp} N\), then \(\operatorname{grade}(M, N) \leq \dim N - \dim M\).
\end{itemize}

    \item By item (1), we always have $\operatorname{depth} R \leq \operatorname{grade} M + \dim M \leq \dim R$.
    \item \cite[Proposition 1.2.10 (a) and (e)]{BH} The following equalities hold:  
    \[
\begin{array}{lll}
\text{grade}(M, N) & =\inf \{ \operatorname{depth} \, N_{\mathfrak{p}}   :  \mathfrak{p} \in \text{Supp} \, M \} \\ &=\inf \{ \operatorname{depth} \, N_{\mathfrak{p}}   :  \mathfrak{p} \in \text{Supp} \, M \cap \text{Supp} \, N \}.
\end{array}
\]
\end{enumerate}
\end{facts}

\begin{newclaim}

For a complex $$X: \cdots \stackrel{\partial_{i+2}}{\longrightarrow} X_{i+1} \stackrel{\partial_{i+1}}{\longrightarrow  } X_i \stackrel{\partial_{i}}{\longrightarrow
      } X_{i-1}   \longrightarrow  \cdots $$ of $R$-modules, we set for each integer $i$, $\operatorname{Z}_i(X)=\ker \partial_i$ and $\operatorname{B}_i(X)= \operatorname{Im} \partial_{i+1}$ and $\operatorname{H}_i(X)=\operatorname{Z}_i(X)/\operatorname{B}_i(X)$.  Moreover, we set: \begin{align*}
&\left\{ 
    \begin{aligned}
        \sup X &= \sup \{ i \in \mathbb{Z} : X_i \neq 0 \},\\
        \inf X &= \inf \{ i \in \mathbb{Z} : X_i \neq 0 \},
    \end{aligned}
\right. \quad 
\left\{
    \begin{aligned}
        \text{hsup } X &= \sup \{ i \in \mathbb{Z} : \operatorname{H}_i(X) \neq 0 \},\\
        \text{hinf } X &= \inf \{ i \in \mathbb{Z} : \operatorname{H}_i(X) \neq 0 \}.
    \end{aligned}
\right.
\end{align*}
The \textit{length} of $X$ is defined to be  $\operatorname{length} X= \sup X - \inf X$. We say that $X$ is \textit{bounded}, if $\operatorname{length} X < \infty$. We say that $X$ is \textit{bounded below} if $\inf X > -\infty$ and $X$ is \textit{bounded above} if $\sup X < \infty$.
\end{newclaim}

\begin{newclaim} Let $M$ be an $R$-module.  We denote by  $\pd_R M$ (resp.  $\operatorname{id}_R M$) the projective (resp. injective) dimension of $M$. 
\end{newclaim}

\subsection{Quasi-homological dimensions.} The definition of quasi-projective dimension (resp. quasi-injective dimension) was introduced by Gheibi, Jorgensen and Takahashi \cite{GJT} (resp. Gheibi \cite{Gheibi})  as a generalization of the classical notion of projective dimension (resp. injective dimension).

\begin{definition}
    Let $M$ be an $R$-module.
    \begin{enumerate}
        \item A \textit{quasi-projective resolution} of $M$ is a bounded below complex $P$ of projective $R$-modules such that for all $i \geq \inf P$ there exist non-negative integers $a_i$, not all zero, such that $H_i(P) \cong M^{\oplus a_i}$. We define the \textit{quasi-projective dimension of $M$} by
        \begin{align*}
            \operatorname{qpd}_R M = \inf \lbrace \sup P - \operatorname{hsup} P : P \text{ is a bounded quasi-projective resolution of } M \},
        \end{align*}
        if  $M\not=0$, and $\operatorname{qpd}_R M=-\infty$ if $M=0$.
        \item A \textit{quasi-injective resolution} of $M$ is a bounded above complex $I$ of injective $R$-modules such that for all $i \leq \sup I$ there exist non-negative integers $b_i$, not all zero, such that $H_i(I) \cong M^{\oplus b_i}$. We define the \textit{quasi-injective dimension} of $M$ by 
  \begin{align*}
            \operatorname{qid}_R M = \inf \lbrace \operatorname{hinf} I - \operatorname{inf} I :  I \text{ is a bounded quasi-injective resolution of } M \},
        \end{align*}
        if  $M\not=0$, and $\operatorname{qid}_R M=-\infty$ if $M=0$.
    \end{enumerate}
    One has $\text{qpd}_R M = \infty$ (resp. $\operatorname{qid}_R M=\infty$) if and only if $M$ does not admit a bounded quasi-projective resolution (resp. quasi-injective resolution). We remark that $\operatorname{qpd}_R M$ is finer than projective dimension, in the sense that there is always an inequality $\operatorname{qpd}_R M \leq \operatorname{pd}_R M$ and equality holds when $\operatorname{pd}_R M$ is finite (see \cite[Corollary 4.10]{GJT}). The same works for the quasi-injective dimension, that is, $\operatorname{qid}_R M \leq \operatorname{id}_R M$ and equality holds when $\operatorname{id}_R M$ is finite, under the assumption that $R$ is local (see \cite[Corollary 3.3]{Gheibi}).
\end{definition}
We observe that results concerning the quasi-projective and quasi-injective dimensions can be found in \cite{GJT, Gheibi}. Below, we state some very useful results related to the quasi-homological dimensions to be used in this paper.
\begin{newclaim}\label{prop41}
Let $R$ be a local ring. A complex $(X,\partial)$ of free $R$-modules of finite rank is called \textit{minimal} if $\partial_i \otimes_R k=0$ for all $i$, where $k$ is the residue field of $R$. If $M$ is a non-zero $R$-module with $\operatorname{qpd}_R M<\infty,$ then  there exists a finite minimal quasi-projective resolution $F$ of $M$ such that $\mathrm{qpd}_R M = \sup F - \mathrm{hsup} F$ (see \cite[Proposition 4.1]{GJT}).
\end{newclaim}

\begin{theorem}\cite[Theorem 4.4]{GJT}\label{rem:ABF} 
Let $R$ be a local ring, and let $M$ be an $R$-module of finite quasi-projective dimension. Then
\[
\operatorname{qpd}_R M = \operatorname{depth} R - \operatorname{depth} M.
\]
In particular, if $M\not=0$, then $\operatorname{depth} M\leq \operatorname{depth} R$ and $\operatorname{qpd}_R M\leq \operatorname{depth} R.$
\end{theorem}
\begin{theorem}\cite[Theorem 4.11]{GJT}\label{depthformula} 
Let $R$ be a local ring, and let $M$ and $N$ be $R$-modules. Suppose that $M$ has finite quasi-projective dimension and $\operatorname{Tor}_{i}^R(M,N)=0$ for all $i>0$. Then 
\[
\operatorname{depth} M + \operatorname{depth} N = \operatorname{depth} R + \operatorname{depth} (M \otimes_R N).
\]
\end{theorem}
\begin{theorem}\cite[Theorem 3.2]{Gheibi}\label{inj} Let $R$ be a local ring and let $M$ be a non-zero $R$-module. If $\operatorname{qid}_R M < \infty$, then $\operatorname{qid}_R M = \operatorname{depth} R$.
    
\end{theorem}
\subsection{Gorenstein dimension} The notion of Gorenstein dimension was introduced by Auslander \cite{Auslander1967} and developed by Auslander and Bridger in \cite{AuBr}. For an $R$-module $M$, set $M^\ast=\Hom_R(M,R)$. 
\begin{definition}
Let $M$ be an $R$-module.
\begin{enumerate}
    \item We say that $M$ is \textit{totally reflexive} if $M$ is reflexive and $\operatorname{Ext}^i_R(M, R) = 0 = \operatorname{Ext}^i_R(M^*, R)$ for all $i > 0$.
    \item The \textit{Gorenstein dimension} of $M$, denoted by $\operatorname{G-dim}_R M$, is defined to be the infimum of all non-negative integers $k$ such that there exists an exact sequence
    $$0 \to G_k \to \cdots \to G_0 \to M \to 0$$
    where each $G_i$ is totally reflexive.
\end{enumerate}
\end{definition}
We can observe that $\operatorname{G-dim}_R M =0$ if and only if $M$ is totally reflexive. 

\subsection{Complete intersection dimension} The notion of complete intersection dimension was introduced by Avramov, Gasharov and Peeva \cite{avramov2}. 
\begin{definition} Let $R$ be a local ring. A diagram of local ring maps $R \to R' \twoheadleftarrow S$ is called a \textit{quasi-deformation} if $R \to R'$ is flat and the kernel of the surjection $R' \twoheadleftarrow S$ is generated by a regular sequence on $S$. The \textit{complete intersection dimension} of an $R$-module $M$ is defined as follows: 
\begin{align*}
\cdim_R M = \inf \lbrace \pd_S (M \otimes_R R')- \operatorname{pd}_S R' : R \rightarrow R' \twoheadleftarrow S \text{ is a quasi-deformation}\rbrace.
\end{align*}
\end{definition}

    




 \section{Quasi-projective dimension and Ischebeck's formula}\label{section3}

In \cite[Theorems 6.2(4) and 6.12]{GJT}, the authors studied the vanishing of the $\operatorname{Ext}$-modules $\operatorname{Ext}_R^i(M,N)$ in two scenarios: when \( M \) has finite quasi-projective dimension, and when \( M \) has finite Gorenstein dimension while \( N \) has finite quasi-projective dimension. Inspired by these results and \cite[Theorem 4.2]{ArayaYoshino}, we investigate the value of \( \operatorname{P}_R(M,N) \) in these cases, assuming that \( \operatorname{P}_R(M,N) < \infty \) and that $R$ is local.  More precisely, the goal of this section is to prove that for non-zero $R$-modules $M$ and $N$ such that $\operatorname{P}_R(M, N)<\infty$ the equality 
\begin{equation}\label{Isck}
    \operatorname{P}_R(M,N)=\depth R-\depth M 
\end{equation}
holds in each one of the following cases:
\begin{enumerate}
\item $M$ has finite quasi-projective dimension. 
\item  $M$ has finite Gorenstein dimension and $N$ has finite quasi-projective dimension.
      \end{enumerate}

\subsection{The first case} In this subsection, we prove the equality \eqref{Isck} considering that $\operatorname{P}_R(M,N)<\infty$ holds in the case (1). The following lemma will be needed. 

\begin{lemma}\label{lemmavanishing}
    Let $R$ be a local ring, and let $M$ and $N$ be non-zero $R$-modules such that $\operatorname{P}_R(M,N)<\infty$. One then has
the inequality $\operatorname{P}_R(M,N) \leq \operatorname{qpd}_R M$.
\end{lemma}
\begin{proof}
We may assume that $\operatorname{qpd}_R M < \infty$. Then, by \cite[Corollary 6.4]{GJT}, we have that $\operatorname{Ext}_R^i(M,N)=0$ for all $i \geq \operatorname{qpd}_R M +1$ and consequently $\operatorname{P}_R(M,N) \leq \operatorname{qpd}_R M$.
\end{proof}
\begin{theorem}\label{teo:qproj}
    Let $R$ be a local ring, and let $M$ and $N$ be non-zero $R$-modules such that $\operatorname{P}_R(M,N)<\infty$. If $\operatorname{qpd}_R M<\infty$, then $$\operatorname{P}_R(M,N)=\operatorname{qpd}_R M=\operatorname{depth} R - \operatorname{depth} M.$$ 
\end{theorem}
\begin{proof}
The second equality is just the Auslander-Buchsbaum formula for quasi-projective dimension (Theorem \ref{rem:ABF}). Since $M$ has finite quasi-projective dimension, by  \ref{prop41}, there exists a minimal quasi-projective resolution 
$$F: \,\, 0 \rightarrow F_s  \rightarrow \cdots \rightarrow F_h \xrightarrow{\partial_h} \cdots $$
of $M$, where $s=\sup F,$ $h=\operatorname{hsup} F$ and $r:=\operatorname{qpd}_R M=s-h$. Letting $C=\operatorname{Coker}(\partial_{h+1})$, we obtain an exact sequence 
$$0 \rightarrow F_s \xrightarrow{\partial_s} \cdots \rightarrow F_h \rightarrow C \rightarrow 0,$$
which is a minimal free resolution of $C$, so that $\operatorname{pd}_R C=s-h=:r$. 

For all $i \in \mathbb{Z}$, setting $Z_i=Z_i(F)$ and $B_i=B_i(F)$, note that there are exact sequences
\begin{align*}
0 \rightarrow Z_i \rightarrow F_i \rightarrow B_{i-1} \rightarrow 0 \text{ and } 0 \rightarrow B_i \rightarrow Z_i \rightarrow H_i(F) \rightarrow 0.
\end{align*}
Then we see that
\[
\operatorname{Ext}^j_R(B_i, N) \cong \operatorname{Ext}^j_R(Z_i, N) \cong \operatorname{Ext}^{j+1}_R(B_{i-1}, N) \cong \cdots \cong \operatorname{Ext}^{j+i+1-\inf F}_R(B_{\inf F - 1}, N) = 0
\]
for all \( i \geq \inf F \) and \( j > p:=\operatorname{P}_R(M,N) \), where the last equality is due to \( B_{\inf F-1} = 0 \). Thus, $\operatorname{Ext}^j_R(B_i, N)=0$ for all $i$ and $j>p$.
Since $h=\operatorname{hsup} F$, then there exists a short exact sequence
\begin{equation*} 
    0 \longrightarrow M^{\oplus a_h} \longrightarrow C \longrightarrow \operatorname{B}_{h-1} \longrightarrow 0
\end{equation*}
for a positive integer $a_h$. Therefore, as $\operatorname{Ext}^j_R(B_{h-1}, N)=0=\Ext_R^j(M,N)$ for all $j>p$, then $\Ext_R^j(C,N)=0$ for all $j>p$. Thus, if $r > p$, then $\Ext_R^r(C,N)=0$, which  contradicts \cite[p. 154, Lemma 1(iii)]{Matsu} as $\operatorname{pd}_R C=r$. Therefore $r\leq p$, and it follows from Lemma \ref{lemmavanishing} that $r=p$. 
\end{proof}
Araya and Yoshino (\cite[Theorem 4.2]{ArayaYoshino}) proved the version of the above theorem with complete intersection dimension instead of quasi-projective dimension. As an application of Theorem \ref{teo:qproj}, we recover this result. We mention that for a local ring $(R,\mathfrak{m},k)$, one always has that $\operatorname{qpd}_R k < \infty$ (see \cite[Proposition 3.6(1)]{GJT}). However, $\operatorname{CI-dim}_R k< \infty$ if and only if $R$ is a complete intersection ring.
\begin{corollary}\cite[Theorem 4.2]{ArayaYoshino} Let $R$ be a local ring, and let $M$ and $N$ be non-zero $R$-modules such that $\operatorname{P}_R(M,N)< \infty$. If $\operatorname{CI-dim}_R M < \infty$, then 
\begin{align*}
\operatorname{P}_R(M,N)=\operatorname{depth} R -\operatorname{depth} M. 
\end{align*}
\end{corollary}
\begin{proof}
Since $\cdim_R M < \infty$, there exists a quasi-deformation $R \rightarrow R' \twoheadleftarrow S$ such that $\operatorname{pd}_S (M \otimes_R R')< \infty$. Set $(-)^\prime=(-) \otimes_R R'$. By the flatness of $R \to R'$, we have $\operatorname{P}_R(M,N)=\operatorname{P}_{R'}(M',N')$ and $\operatorname{depth} R - \operatorname{depth} M = \operatorname{depth} R' - \operatorname{depth} M'$. Thus, we may assume that $R=R'$, i.e., $R = S/(\boldsymbol{x})$, where $\boldsymbol{x}$ is a regular sequence on $S$. Since $\operatorname{pd}_S M < \infty$, then $\operatorname{qpd}_R M < \infty$, by \cite[Proposition 3.7]{GJT}. Thus, the desired equality follows directly by Theorem \ref{teo:qproj}.
\end{proof}

The New Intersection Theorem (that was proved by Peskine and Szpiro \cite{PS} for some cases, and by Hochster \cite{Hoch} in the equicharacteristic case and finally by Roberts \cite{Roberts} in the mixed characteristic case), shows that if \( M \) and \( N \) are \( R \)-modules, one has the inequality \( \dim N \leq \operatorname{pd}_R M + \dim (M \otimes_R N) \). However, the version with quasi-projective dimension of this inequality, generally does not hold (see \cite[Example 4.7]{arxiv}). The following proposition provides, as a corollary, a context in which this is true.

\begin{proposition}\label{prop:int}
Let $R$ be a local ring, and let $M$ and $N$ be non-zero $R$-modules such that $\operatorname{P}_R(M,N)< \infty$. Then we have the inequality \[
\dim \operatorname{Ext}^i_R(M, N) + i \leq \operatorname{qpd}_R M + \dim(M \otimes_R N) \quad \text{for all} \ i.
\]
\end{proposition}
\begin{proof} 
We may assume $\operatorname{qpd}_RM < \infty$. Since 
\[
\operatorname{Supp} \operatorname{Ext}^i_R(M, N) \subseteq \operatorname{Supp} M \cap \operatorname{Supp} N = \operatorname{Supp}(M \otimes_R N) \quad \text{for all} \ i,
\]
we have 
\begin{equation*}
\dim \operatorname{Ext}^i_R(M, N) \leq \dim(M \otimes_R N) \quad \text{for all} \ i.
\end{equation*}

Now, since $\operatorname{qpd}_R M < \infty$, for integers $i$ such that $\operatorname{Ext}_R^i(M,N) \neq 0 $ note that $i \leq \operatorname{qpd}_R M$ by Theorem \ref{teo:qproj}, and then we obtain the desired inequality. For those $i$ for which $\Ext_R^i(M,N)=0$, we have that $\dim \Ext_R^i(M,N)=-\infty$ (see \ref{convetion}), so the inequality is trivially satisfied.
\end{proof}



\begin{corollary}\label{intersection}
Let $R$ be a local ring, and let $M$ and $N$ be non-zero $R$-modules such that $\operatorname{P}_R(M,N)< \infty$. Then we have the inequality  
\begin{align*}
\depth N \leq \operatorname{qpd}_R M + \dim (M \otimes_R N).
\end{align*}
In particular, if $N$ is Cohen-Macaulay, then $\dim N \leq \operatorname{qpd}_R M + \dim (M \otimes_R N)$.
\end{corollary}
\begin{proof}
Since $\operatorname{P}_R(M,N) < \infty$, it follows from \cite[Proposition 3.5]{Foxby1979} and \cite[Lemma 2.16] {Yassemi1} \footnote{In the statement of \cite[Lemma 2.16]{Yassemi1}, the modules $M$ and $N$ are required to be non-zero, although this is not stated explicitly.} that 
$\depth N \leq \dim \operatorname{Ext}_R ^i(M,N)+i$ for some $i$. Hence, the desired inequality follows from Proposition \ref{prop:int}.
\end{proof}

 \subsection{The second case.} Next, we prove that $\operatorname{P}_R(M,N) = \depth M - \depth R$, assuming that $\operatorname{P}_R(M,N) < \infty$ holds in case (2). To this end, we first establish the following lemma.

 \begin{lemma}\label{lemm2}
Let $R$ be a local ring, and let $M$ and $N$ be non-zero $R$-modules such that $\operatorname{P}_R(M,N)< \infty$ and $\operatorname{qpd}_R N < \infty$. One then has
the inequality $\operatorname{P}_R(M,N)\leq \operatorname{G-dim}_RM$. 
 \end{lemma}
 \begin{proof}
We may assume that $\operatorname{G-dim}_R M < \infty$. By contradiction, assume that $p:=\operatorname{P}_R(M,N)>\operatorname{G-dim}_R M$. Since $\operatorname{Ext}_R^i(M,N)=0$ for all $i>p$, then \cite[Lemma 6.11(2)]{GJT} implies that $\operatorname{Ext}_R^p(M,N)=0$. This is a contradiction to the definition of $p$. 
 \end{proof}

\begin{fact}{\cite[Remark 3.3]{ChristensenIyengar2007}}\label{fact9}
Let $M$ be an $R$-module with $\operatorname{G-dim}_RM < \infty$. Then, there exists a short exact sequence $ 0 \to M \to H \to X \to 0$, with $X$ totally reflexive and $\operatorname{pd}_R H = \operatorname{G-dim}_R M$.  
\end{fact}

\begin{theorem}\label{g-dim}
     Let $R$ be a local ring, and let $M$ and $N$ be non-zero $R$-modules such that $\operatorname{P}_R(M,N)<\infty$. If $\operatorname{G-dim}_R M<\infty$ and $\operatorname{qpd}_R N<\infty$, then $$\operatorname{P}_R(M,N)=\operatorname{G-dim}_R M= \operatorname{depth} R -\operatorname{depth} M.$$
\end{theorem}
 \begin{proof}
Let $p=\operatorname{P}_R(M,N)$. The second equality is just the classical Auslander-Bridger formula (\cite[Theorem 29]{masiek}). We will prove the first equality by induction on $p$. Suppose $p=0$, that is $\Ext_R^i(M,N)=0$ for all $i>0$. Since $\operatorname{G-dim}_R M< \infty$, we can use a Gorenstein dimension approximation of $M$ (see Fact \ref{fact9}), that is, a short exact sequence \begin{equation}\label{q1zp}
    0 \to M \to H \to X \to 0,
\end{equation}
where $X$ is totally reflexive and $\operatorname{pd}_R H=\operatorname{G-dim}_R M$. Since  $\operatorname{pd}_R H<\infty$, then $\operatorname{P}_R(H,N)<\infty$. Hence, using the long sequence $\Ext(-, N)$  induced from \eqref{q1zp}, it follows that $\operatorname{P}_R(X, N)<\infty$.  Thus, applying  Lemma \ref{lemm2}, we have that $\Ext_R^i(X,N)=0$ for all $i>0$. Then using again the long sequence $\Ext(-, N)$  induced from \eqref{q1zp}, it follows that $\Ext_R^i(H,N)=0$ for all $i>0$. Therefore, by \cite[p. 154, Lemma 1(iii)]{Matsu},  $H$ is free. Consequently, the equality $\operatorname{pd}_R H=\operatorname{G-dim}_R M$ says that $\operatorname{G-dim}_R M=0=p$.  

 Now, suppose that $p>0$. We see that  $\operatorname{P}_R(\Omega^1M,N)=p-1$ and, by Lemma \ref{lemm2}, $M$ is not totally reflexive. We also have that $\operatorname{G-dim}_R( \Omega^1 M)=\operatorname{G-dim}_R M -1$. By induction, we obtain that $\operatorname{P}_R(\Omega^1M, N)=\operatorname{G-dim}_R(\Omega^1 M)$. Thus, $p-1=\operatorname{G-dim}_R M-1,$ whence $p=\operatorname{G-dim}_R M$, as desired.
\end{proof}
As an application, we derive the following proposition, which generalizes \cite[Theorem 1.4]{Yassemi1}.


\begin{proposition}
Let $R$ be a local ring, and let $M$ and $N$ be non-zero $R$-modules such that $\operatorname{P}_R(M,N)<\infty$ and $\operatorname{qpd}_R N < \infty$. One then has
the inequality
 \[
\dim \operatorname{Ext}^i_R(M, N) + i \leq \operatorname{G-dim}_R M + \dim(M \otimes_R N) \quad \text{for all} \ i.
\]
\end{proposition}
\begin{proof}
The proof follows the same lines as the proof of Proposition \ref{prop:int}, substituting \(\operatorname{G-dim}_R M\) with \(\operatorname{qpd}_R M\) and Theorem \ref{g-dim} with Theorem \ref{teo:qproj}.
\end{proof}


\section{A generalization of a result by Ischebeck}\label{section4}

Let $R$ be a local ring, and let $M$ and $N$ be non-zero $R$-modules. Ischebeck \cite[2.6]{isch} proved that \( \operatorname{P}_R(M, N) = \operatorname{depth} R - \operatorname{depth} M \), assuming that \( N \) has finite injective dimension. A similar result was established by Sazeedeh \cite{saz}, when \( M \) has finite injective dimension and \( N \) has finite Gorenstein injective dimension. In the main result of this section (Theorem \ref{isch}), we generalize Ischebeck's  result considering quasi-injective dimension in place of injective dimension and adding the condition that $\operatorname{P}_R(M,N) < \infty$. 

\begin{lemma}\label{lemmadepthzero}
Let $R$ be a local ring, and let $M$ and $N$ be non-zero $R$-modules such that $\operatorname{P}_R(M,N) < \infty$. If $\operatorname{qid}_R N < \infty$, then $\operatorname{depth} M \leq \operatorname{depth} R$.
\end{lemma}
\begin{proof}
Let $\boldsymbol{x}$ be a maximal $M$-regular sequence. By induction on the length of $\boldsymbol{x}$ and using the long exact sequence $\operatorname{Ext}_R(-,N)$, one can see that $\operatorname{P}_R(M/\textbf{\textit{x}}M,N)<\infty$ and $\operatorname{P}_R(M/\textbf{\textit{x}}M,N)=\operatorname{P}_R(M,N)+\operatorname{depth} M$. 
By \cite[Corollary 3.5]{Gheibi}, we have $\operatorname{P}_R(M/\textbf{\textit{x}}M,N) \leq \operatorname{depth} R$. Combining this with the previous equality, we see that $\operatorname{P}_R(M,N)+\depth M \leq \depth R$. Thus, since $\operatorname{P}_R(M,N)\geq 0$ (see \ref{claim2.3}), we  obtain that   $\operatorname{depth} M \leq \operatorname{depth} R$.
\end{proof}

\begin{lemma}\label{lema:notation}
Let $(R,\mathfrak{m},k)$ be a local ring, and let $M$ and $N$ be non-zero $R$-modules with $\operatorname{P}_R(M,N)< \infty$ and $\operatorname{qid}_R N < \infty$. Let 
$$I: \,\, 0 \rightarrow I_0 \xrightarrow{\partial_0} I_{-1}  \xrightarrow{\partial_{-1}} I_{-2} \rightarrow \cdots  $$
be a bounded quasi-injective resolution of $N$ such that $\operatorname{qid}_R N = \operatorname{hinf} I - \operatorname{inf} I$ and  $\operatorname{sup} I= \operatorname{hsup} I=0$ (such $I$ always exists, by \cite[Remark 2.3(3)]{Gheibi}). Set $s=\operatorname{hinf} I$, $Z_i=Z_i(I)$ and $B_i=B_i(I)$, for all $i \in \mathbb{Z}$. Then
\begin{enumerate}
    \item $\operatorname{P}_R(M,N)=\operatorname{P}_R(M,Z_s)$.
    \item If $t=\operatorname{depth} R >0$, then $\operatorname{Ext}_R^t(k,Z_s) \neq 0$.
\end{enumerate}
\end{lemma}
\begin{proof} (1) Set $p=\operatorname{P}_R(M,N)$ and $p_s=\operatorname{P}_R(M,Z_s)$. For all $i \in \mathbb{Z}$, there are exact sequences
\begin{align}\label{seqs}
0 \rightarrow Z_i \rightarrow I_i \rightarrow B_{i-1} \rightarrow 0 \text{ and } 0 \rightarrow B_i \rightarrow Z_i \rightarrow H_i(I) \rightarrow 0. 
\end{align}
We claim that $\Ext_R^{i>p}(M,B_{-j})=0$ and $\Ext^{i>p}_R(M,Z_{-j})=0$ for all $j \geq 0$. Indeed, it is clear that $\Ext_R^{i > p} (M,B_0)=0$ since $B_0=0$. By using the exact sequence $0 \rightarrow Z_0 \rightarrow I_0 \rightarrow B_{-1} \rightarrow 0 $ and noting that $Z_0 \cong  N^{\oplus b_0}$ for some positive integer $b_0$, the fact that $\Ext^{i>p}_R (M,N)=0$ allow us to obtain that $\Ext^{i>p}_R(M,B_{-1})=0$. Now, considering the exact sequence $0 \rightarrow B_{-1} \rightarrow Z_{-1} \rightarrow H_{-1}(I) \rightarrow 0$, we can see that $\Ext_R^{i>p}(M,Z_{-1}) =0$. Therefore, using the exact sequences of \eqref{seqs},  we conclude inductively the vanishing of desired Ext-modules. In particular, we have $\operatorname{Ext}_R^{i>p}(M,Z_s)=0$ and thus $p_s \leq p$. To prove the desired equality,  suppose by contradiction that $p_s < p$.  Then  $\operatorname{Ext}_R^p(M,Z_s)=0$. Since $\operatorname{H}_s(I)\cong 
 N^{\oplus b_s}$ for some positive integer $b_s$, there exists a short exact sequence
$$0 \rightarrow B_s \rightarrow Z_s \rightarrow  N^{\oplus b_s} \rightarrow 0.$$
Thus from its long exact sequence $\operatorname{Ext}_R(M,-)$, using the fact that $\operatorname{Ext}_R^{p+1}(M,B_s)=0$, we see that $\Ext_R^p(M,N)=0$, a contradiction to definition of $p$. Hence, $p=p_s$.

(2) We have that $t=\operatorname{qid}_R N = \operatorname{id}_R Z_s$, by Theorem \ref{inj}. For a prime ideal $\mathfrak{p} \neq \mathfrak{m}$ of $R$ and $x \in \mathfrak{m} \backslash \mathfrak{p}$, consider the exact sequence $0 \longrightarrow R / \mathfrak{p} \xrightarrow{x} R / \mathfrak{p}$. This yields an exact sequence
$$
\operatorname{Ext}_R^t\left(R / \mathfrak{p}, Z_s\right) \xrightarrow{x} \operatorname{Ext}_R^t\left(R / \mathfrak{p}, Z_s\right) \longrightarrow 0 .
$$
Since $t>0$, then \cite[Lemma 3.1]{Gheibi} says that  $\operatorname{Ext}_R^t(R/\mathfrak{p},Z_s)$ is finitely generated, and therefore, by Nakayama's Lemma, we have that $\operatorname{Ext}_R^t\left(R / \mathfrak{p}, Z_s\right)=0$. Since $\operatorname{id}_R Z_s=t$, we must then have $\operatorname{Ext}_R^t\left(k, Z_s\right) \neq 0$ by \cite[Corollary 3.1.12]{BH}. 
\end{proof}
Note that, in fact, item (2) of the lemma above is deduced from the proof \cite[Theorem 3.2]{Gheibi}.


\begin{theorem}\label{isch}
Let $(R,\mathfrak{m},k)$ be a local ring, and let $M$ and $N$ be non-zero $R$-modules with $\operatorname{P}_R(M,N)<\infty$. If $\operatorname{qid}_R N < \infty$, then 
\begin{align*}
    \operatorname{P}_R(M,N)=\operatorname{depth} R- \operatorname{depth} M.
\end{align*}
\end{theorem}

\begin{proof}
If $\operatorname{depth} R =0$, Lemma \ref{lemmadepthzero} yields that  $\operatorname{depth} M=0$ and  $\operatorname{P}_R(M,N)=0$ by \cite[Corollary 3.5]{Gheibi}. Thus, the desired equality holds. So, we may assume  $t=\operatorname{depth} R >0$. Consider the notation introduced in Lemma \ref{lema:notation}. By Theorem \ref{inj}, we have that $t=\operatorname{qid}_R N = \operatorname{id}_RZ_s$. By Lemma \ref{lema:notation}(1), it is enough to prove that $$\operatorname{P}_R(M,Z_s)=  \operatorname{depth} R - \operatorname{depth} M.$$ 
We prove it by induction on $\depth M$. If $\operatorname{depth} M=0$, then there exists an exact sequence $$0 \rightarrow k \rightarrow M \rightarrow C \rightarrow 0.$$
    This induces an exact sequence $$\Ext_R^t(M, Z_s) \rightarrow \Ext_R^t(k, Z_s) \rightarrow \Ext_R^{t+1}(C,Z_s)=0.$$
    Then $\Ext_R^t(M,Z_s)\not=0$ since $\operatorname{Ext}_R^t(k,Z_s) \neq 0$ by Lemma \ref{lema:notation}(2). Additionally, we have that $\Ext_R^i(M, Z_s)=0$ for $i>t$ as $t=\operatorname{id}_R Z_s$. Thus, $\operatorname{P}_R(M,Z_s)=t=\depth R$, as desired.  

      Now, assume that $r:=\depth M>0$. Let $x \in \mathfrak{m}$ be an $M$-regular element. Consider the exact sequence 
      \begin{align}\label{seq5}
          0 \rightarrow M \xrightarrow{x} M \rightarrow M/xM \rightarrow 0. 
          \end{align}
        Since $\operatorname{P}_R(M,N)<\infty$, using its long exact sequence $\Ext_R(-,N)$, we see that \linebreak $\operatorname{P}_R(M/xM, N)<\infty$. Then, by induction, $\operatorname{P}_R(M/xM,Z_s)=t-r+1$. The exact sequence \eqref{seq5}
induces for all $i\geq 0$ an exact sequence 
\begin{equation}\label{p9q8}
        \Ext_R^i(M, Z_s) \xrightarrow{x} \Ext_R^i(M, Z_s) \rightarrow \Ext_R^{i+1}(M/xM,Z_s) \rightarrow \operatorname{Ext}_R^{i+1} (M,Z_s)
    \end{equation}
    which consists of finitely generated $R$-modules when $i>0$ by \cite[Lemma 3.1]{Gheibi}. By Lemma \ref{lemmadepthzero}, we have that $t-r \geq 0 $. If we consider $i>t-r$ in \eqref{p9q8}, then $\operatorname{Ext}_R^{i+1}(M/xM,Z_s)=0$ as $\operatorname{P}_R(M/xM,Z_s)=t-r+1$, whence, by Nakayama's Lemma,  $\Ext_R^{i}(M, Z_s)=0$.  If we take $i=t-r$ in \eqref{p9q8}, then $\operatorname{Ext}_R^{t-r}(M,Z_s)\not=0$ since  $\Ext_R^{t-r+1}(M/xM,Z_s)\not=0=\Ext_R^{t-r+1}(M,Z_s)$. Therefore, $\operatorname{P}_R(M,Z_s)=t-r$.

\end{proof}
The following example of Khatami and Yassemi \cite{Khatami1} demonstrates that the condition \( \operatorname{P}_R(M, N) < \infty \) in Theorems \ref{teo:qproj}, \ref{g-dim} and \ref{isch} cannot be omitted.
\begin{example}
Let $(R,\mathfrak{m},k)$ be a Gorenstein local ring which is not regular. The residue field $k$ has finite quasi-projective, quasi-injective and Gorenstein dimensions (see \cite[Proposition 3.6(1)]{GJT}, \cite[Proposition 2.8(1)]{Gheibi} and \cite[Theorem 17]{masiek}, respectively). On the other hand, since $R$ is non-regular, we must have that $\operatorname{pd}_R k = \infty$ and $\operatorname{P}_R(k,k)=\infty$.
\end{example}

\section{Grade and quasi-homological dimensions}\label{section5}
In this section, we present results concerning the relation between the grade and the quasi-homological dimensions.



\subsection{Grade, quasi-projective dimension and quasi-perfect modules} In this subsection, we provide a relation between grade and quasi-projective dimension, and introduce a new class of modules, called the quasi-perfect modules. 

It is well known that for an $R$-module $M$ the following chain of inequalities holds: $\operatorname{grade} M \leq \operatorname{G-dim}_R M \leq \operatorname{pd}_R M$. Then the following question arises naturally. 


\begin{question}
Let $M$ be an $R$-module. Does the inequality $\operatorname{grade} M \leq \operatorname{qpd}_R M$ hold? 
\end{question}
The next theorem says that the answer to this question is affirmative.

\begin{theorem}\label{theo4.5}
Let $M$ be a non-zero $R$-module. One then has
the inequality  $\operatorname{grade} M \leq \operatorname{qpd}_R M$.
\end{theorem}

\begin{proof} For each  $\mathfrak{p} \in \operatorname{Supp} M$, we have that $\operatorname{grade} M \leq \operatorname{grade}_{R_\mathfrak{p}}M_{\mathfrak{p}}$ and, by \cite[Proposition 3.5(1)]{GJT}, $\operatorname{qpd}_{R_{\mathfrak{p}}}M_{\mathfrak{p}} \leq \operatorname{qpd}_R M$. In view of this, we may assume that $R$ is local. We will proceed by induction on $\operatorname{grade} M$.

If $\operatorname{grade} M=0$, it is clear that the inequality holds. Suppose $\operatorname{grade} M > 0$.  Hence there exists an $R$-regular element $x$ such that $xM=0$ by \cite[p. 129, Theorem 16.6]{Matsu}. We may assume that $\operatorname{qpd}_R M<\infty$. Then, by \cite[Proposition 2.11(1)]{Gheibi}, there exists a positive integer $n$ such that $\operatorname{qpd}_{R/(x^n)} M < \infty$. Then we have \begin{align*}
    \operatorname{qpd}_{R/(x^n)} M &= \depth_{R/(x^n)} R/(x^n)-\depth_{R/(x^n)} M\\&=\depth R-1-\depth_R M \\&=\operatorname{qpd}_R M -1,
\end{align*}
where the first and third equalities follow from the Auslander-Buchsbaum formula for quasi-projective dimension, while the second equality is due to \cite[Exercise 1.2.26(b)]{BH}. Moreover since  $\operatorname{grade} M>0$, from \cite[p. 140, Lemma 2(i)]{Matsu}, one can derive $ \operatorname{grade}_{R/(x^n)} M=\operatorname{grade} M -1$. Hence, by induction, $\operatorname{grade}_{R/(x^n)} M \leq \operatorname{qpd}_{R/(x^n)} M$. Consequently, $\operatorname{grade} M \leq \operatorname{qpd}_R M$.
\end{proof}

\begin{definition}
Let $M$ be an $R$-module. We define the {\it Cohen-Macaulay defect} of \( M \), denoted by \( \operatorname{cmd} M\), as the difference \(\dim M - \operatorname{depth} M\).    
\end{definition}
 
It is clear that $M$ is Cohen-Macaulay as an $R$-module if and only if $\operatorname{cmd} M =0$. The next theorem relates the grade of two modules of finite quasi-projective dimension with their quasi-projective dimensions. 
\begin{proposition}\label{Proposition5.6}
Let \( R \) be a local ring, and let \( M \) and \( N \) be non-zero \( R \)-modules of finite quasi-projective dimensions. Then
\[
\qpd_R M -\operatorname{grade}(M, N) \leq \qpd_R N + \operatorname{cmd} M.
\]
\end{proposition}
\begin{proof}
    By Facts \ref{facts3}(1), we have the inequality $ -\operatorname{grade}(M, N) \leq   \dim M-\operatorname{depth} N$. Thus
   $$\depth R-\depth M -\operatorname{grade}(M, N) \leq   \dim M-\operatorname{depth} N+\depth R-\depth M,$$
   and by Auslander-Buchsbaum formula for quasi-projective dimension, we get: 
   \[
\qpd_R M -\operatorname{grade}(M, N) \leq \qpd_R N + \operatorname{cmd} M.
\]
   \end{proof}

Let \( M \) be an \( R \)-module. The following chain of inequalities  
\[
\operatorname{grade} M \leq \operatorname{qpd}_R M \leq \operatorname{pd}_R M
\]
suggests a natural extension of the concept of a perfect module, motivating the introduction of a new notion based on the quasi-projective dimension. 
\begin{definition}
Let \( M \) and \( N \) be \( R \)-modules with \(\operatorname{qpd}_R M < \infty\). We say that \( M \) is \( N \)-\textit{quasi-perfect} if \(\operatorname{qpd}_R M = \operatorname{grade}(M, N)\). In particular, if \( N = R \), then \( M \) is called \emph{quasi-perfect} if \(\operatorname{qpd}_R M = \operatorname{grade} M\). 
\end{definition}

\begin{remark}
    (1) If $(R, \mathfrak{m},k)$ is a local ring, then $k$ is a quasi-perfect module (see \cite[Proposition 3.6]{GJT}). More generally, let $R$ be a ring and let $I$ be an ideal of $R$ such that the Koszul complex with respect to a system of generators of $I$ is a quasi-projective resolution of $R/I$ (e.g., when $I$ is a complete intersection or, more broadly, a quasi-complete intersection), then the $R$-module $R/I$ is always quasi-perfect (see \cite[Theorem 7.4(a)]{GJT}). 
    
    (2) There are quasi-perfect modules that are not perfect neither $G$-perfect. For instance, the residue field $k$ is quasi-perfect but is not perfect (resp. $G$-perfect) unless that $R$ is regular (resp. Gorenstein).


    (3) Let $R$ be a local ring. If $\operatorname{P}_R(M,R)<\infty$ and $M$ is quasi-perfect, then $M$ is $\operatorname{G-}$perfect. As shown using \cite[Proposition 6.14]{GJT}, the classical Auslander-Bridger formula, and the Auslander-Buchsbaum formula for quasi-projective dimension.
\end{remark}

Next, we establish some results for quasi-perfect modules similar to the classic results for perfect modules.

We start presenting a relation between Cohen-Macaulayness and quasi-perfectness. To motivate this, we remind a well-known fact: If $R$ is a Cohen-Macaulay local ring, and $M$ is an $R$-module of finite projective dimension, then $M$ is a Cohen-Macaulay $R$-module if and only if $M$ is perfect (see e.g \cite[Theorem 2.1.5]{BH}). Next proposition generalizes this fact in the context of quasi-perfect $R$-modules.
\begin{proposition}\label{prop:eqv}
Let $R$ be a local ring, and let $M$ be a non-zero $R$-module with $\operatorname{qpd}_R M < \infty$. Then the following statements hold:
\begin{enumerate}
    \item If $M$ is Cohen-Macaulay, then $M$ is quasi-perfect.
    \item If $R$ is Cohen-Macaulay and $M$ is quasi-perfect, then $M$ is Cohen-Macaulay.
\end{enumerate}
\end{proposition}
\begin{proof}
\begin{enumerate}
    \item By Theorem \ref{theo4.5}, we have that $\operatorname{grade} M \leq \operatorname{qpd}_R M$ and, since $M$ is Cohen-Macaulay, we have that $\operatorname{qpd}_R M \leq \operatorname{grade} M$ by Proposition \ref{Proposition5.6}. Thus $\operatorname{grade} M = \operatorname{qpd}_R M $, as desired. 
    \item Since $R$ is Cohen-Macaulay and $M$ is quasi-perfect, we have
    \[
    \begin{array}{lll}
       \dim M &= \operatorname{depth} R - \operatorname{grade} M \quad &\text{(by Facts \ref{facts3} (2))} \\
              &= \operatorname{depth} R - \text{qpd}_R M \quad &\text{(by hypothesis)} \\
              &= \operatorname{depth} M \quad & \text{(by Theorem \ref{rem:ABF})}.
    \end{array}
    \]
\end{enumerate}
\end{proof}

The notion of a quasi-perfect ideal was introduced in \cite{Gheibi} as follows. An ideal $I$ of $R$ is said to be \textit{quasi-perfect} if $\operatorname{grade} I:= \operatorname{grade} R/I=\operatorname{qpd}_R R/I$. Thus, observe that an ideal \( I \) of \( R \) is quasi-perfect in the sense of \cite{Gheibi} if and only if the \( R \)-module \( R/I \) is a quasi-perfect \( R \)-module. The following corollary generalizes \cite[Corollary 7.6]{GJT}.

\begin{corollary}
Let $R$ be a Cohen-Macaulay local ring, and let $I$ be a quasi-perfect ideal of $R$. Then $I$ is a Cohen-Macaulay ideal of $R$.
\end{corollary}

Let $M$ and $N$ be $R$-modules. We set $\operatorname{q}^R(M,N):=\sup \lbrace i \geq 0 : \operatorname{Tor}_i^R(M,N) \neq 0 \rbrace$.

\begin{theorem} \label{theorem1.8}
Let $R$ be a local ring. Let $M$ and $N$ be non-zero $R$-modules such that $N$ Cohen-Macaulay and $\operatorname{qpd}_R M < \infty$. If $\operatorname{q}^R(M,N)=0$, then \( M \otimes_R N \) is Cohen-Macaulay if and only if \( M \) is \( N \)-quasi-perfect.
\end{theorem}

\begin{proof}
Since \(\operatorname{Supp}(M \otimes_R N) = \operatorname{Supp} M \cap \operatorname{Supp} N\), from Facts \ref{facts3}(1) and the fact that \( N \) is Cohen-Macaulay, we have
$$\operatorname{depth} N - \dim (M \otimes_R N) = \operatorname{grade}(M \otimes_R N, N).$$

Adding \(-\operatorname{depth} (M \otimes_R N)\) to the equality, we get
$${\rm cmd}(M \otimes_R N) = \operatorname{depth} N - \operatorname{depth} (M \otimes_R N) - \operatorname{grade}(M \otimes_R N, N).$$

On the other hand,
\[
\begin{array}{lll}
\operatorname{grade}(M \otimes_R N, N) &=& \inf \{ \operatorname{depth} N_{\mathfrak{p}}   :  \mathfrak{p} \in \operatorname{Supp}(M \otimes_R N) \}\\
&=& \inf \{ \operatorname{depth} N_{\mathfrak{p}}   :  \mathfrak{p} \in \operatorname{Supp} M \cap \operatorname{Supp} N \}\\
&=& \operatorname{grade}(M, N).
\end{array}
\]

Therefore, we have
$${\rm cmd}(M \otimes_R N) = \operatorname{depth} N - \operatorname{depth} (M \otimes_R N) - \operatorname{grade}(M, N).$$

Now, by the Auslander-Buchsbaum formula and the depth formula both for quasi-projective dimension (see Theorems \ref{rem:ABF} and \ref{depthformula}), we have
$${\rm cmd}(M \otimes_R N) = \operatorname{qpd}_R M - \operatorname{grade}(M, N).$$

The equality shows that \( M \otimes_R N \) is Cohen-Macaulay if and only if \( M \) is \( N \)-quasi-perfect.
\end{proof}
\begin{proposition}\label{prop:mixed}
Let $M$ be a quasi-perfect $R$-module. For any prime ideal $\mathfrak{p} \in \operatorname{Supp} M$ the following are equivalent:
\begin{enumerate}
    \item $\mathfrak{p} \in \operatorname{Ass} M$.
    \item $\operatorname{depth} R_{\mathfrak{p}}= \operatorname{grade} M$.  
\end{enumerate}
Moreover, for every $\mathfrak{p} \in \operatorname{Ass}(M)$, we have $\operatorname{grade} R/\mathfrak{p} = \operatorname{grade}M$.
\end{proposition}

\begin{proof}
Let $\mathfrak{p} \in \operatorname{Supp} M$. Then   \cite[Proposition 3.5(1)]{GJT} and Theorem \ref{theo4.5} yield the following sequence of inequalities:
\begin{align}\label{in}
\operatorname{grade}M \leq \operatorname{grade}M_{\mathfrak{p}} \leq \operatorname{qpd}_{R_{\mathfrak{p}}}M_{\mathfrak{p}} \leq \operatorname{qpd}_R M.
\end{align}
Since $M$ is quasi-perfect, all the inequalities of (\ref{in}) become equalities. In particular, $\operatorname{qpd}_{R_\mathfrak{p}} M_\mathfrak{p}=\operatorname{grade} M_\mathfrak{p}$. 
On the other hand, by the Auslander-Buchsbaum formula for quasi-projective dimension, we get $\operatorname{qpd}_{R_{\mathfrak{p}}}M_{\mathfrak{p}} + \operatorname{depth}M_{\mathfrak{p}} = \operatorname{depth}R_{\mathfrak{p}}$. Then it follows  from the equalities that $\operatorname{depth} M_{\mathfrak{p}} = 0$ if and only if $\operatorname{grade} M = \operatorname{depth}R_{\mathfrak{p}}$. Therefore, the equivalence holds.

For the second part, assume that $\mathfrak{p} \in \operatorname{Ass}(M)$. Then  $ \operatorname{ann}M \subseteq \mathfrak{p} $ and hence,  $\operatorname{grade}M=\operatorname{depth}( \operatorname{ann} M, R)\leq \operatorname{depth}(\mathfrak{p}, R)=\operatorname{grade} R/\mathfrak{p}. $  For the other inequality, note that according to what was proved above, we have $\operatorname{grade} M = \operatorname{depth} R_{\mathfrak{p}}$.  Since $\operatorname{grade} R/\mathfrak{p}\leq \operatorname{depth} R_\mathfrak{p}$, then $\operatorname{grade} R/\mathfrak{p} \leq \operatorname{grade} M$.
\end{proof}

 It is well-known that a local ring is Cohen-Macaulay  if and only if it admits a non-zero Cohen-Macaulay module with finite projective dimension. However, this statement does not hold true when quasi-projective dimension is considered in place of projective dimension. In fact, the residue field of every local ring satisfies these conditions. The next corollary provides criteria for the base ring \( R \) to be Cohen-Macaulay, based on the existence of a Cohen-Macaulay module with finite quasi-projective dimension with an additional hypo\-thesis.
\begin{corollary}\label{pros1}
Let $R$ be a local ring. The following are equivalent:
\begin{enumerate}
    \item $R$ is Cohen-Macaulay.
    \item There exists a non-zero Cohen-Macaulay $R$-module $M$ with $\operatorname{qpd}_R M < \infty$ and $\dim M= \operatorname{dim} R - \operatorname{grade}M$. 
    \item There exists a non-zero Cohen-Macaulay $R$-module $M$ with $\operatorname{qpd}_R M < \infty$ and such that $\dim R/\mathfrak{p}+\operatorname{depth} R_{\mathfrak{p}}=\dim R$ for some $\mathfrak{p} \in \operatorname{Ass} M$.
\end{enumerate}
\end{corollary}
\begin{proof}
The assertions (1) $\Rightarrow$ (2) and (1) $\Rightarrow$ (3) are trivial, by setting $M=R$.

(2) $\Rightarrow$ (1). Since $M$ is Cohen-Macaulay, then $M$ is quasi-perfect, by Proposition \ref{prop:eqv}(1), that is, $\operatorname{grade}M =\operatorname{qpd}_R M$. Therefore, using the Auslander-Buchsbaum formula for quasi-projective dimension and assumption, we see that $\operatorname{dim} R = \operatorname{depth} R$, that is $R$ is Cohen-Macaulay.   

(3) $\Rightarrow$ (2). As before,  since $M$ is Cohen-Macaulay, then $M$ is quasi-perfect. Proposition \ref{prop:mixed} and  \cite[Theorem 2.1.2(a)]{BH} give $\operatorname{depth} R_{\mathfrak{p}}= \operatorname{grade} M$ and $\operatorname{dim} M= \dim R/\mathfrak{p}$, respectively. Substituting these in the equality given by assumption, we obtain $\dim M = \dim R - \operatorname{grade} M$. 

\end{proof}

\subsection{Grade and quasi-injective dimension}
In the main theorem of this subsection, we prove a formula for the grade of finitely generated modules of finite quasi-injective dimension over a local ring. This theorem improves \cite[Theorem 3.7]{Gheibi}.

 
\begin{theorem}\label{teo:forinj}
    Let $R$ be a local ring and let $M$ be a non-zero $R$-module with $\operatorname{qid}_R M <\infty$. Then $$\operatorname{dim} M=\depth R-\operatorname{grade} M.$$
\end{theorem}
\begin{proof}
    We proceed by induction on $\operatorname{grade} M.$ If $\operatorname{grade} M=0$, from Facts \ref{facts3}(2), we have that $\operatorname{depth}R\leq \operatorname{dim}M$. The opposite inequality is also valid by \cite[Theorem 3.7]{Gheibi}. Thus, we have \(\operatorname{depth} R = \operatorname{dim} M\), and the desired equality holds.

    Now, assume that $\operatorname{grade} M>0$. Then, by \cite[p. 129, Theorem 16.6]{Matsu}, there exists an $R$-regular element $x$ with $xM=0$.  Since $\operatorname{qid}_R M<\infty$, then  $\operatorname{qid}_{R/(x^n)} M<\infty$ for some $n>0$ by \cite[Proposition 2.11(2)]{Gheibi}. As $\operatorname{grade}(M)>0$, we can observe from \cite[p. 140, Lemma 2(ii)]{Matsu} that $\operatorname{grade}_{R/(x^n)} M=\operatorname{grade} M-1$. Therefore, by induction, we have:  $$\operatorname{dim}_{R/(x^n)} M=\depth_{R/(x^n)} R/(x^n)-\operatorname{grade}_{R/(x^n)} M. $$
    Thus, 
    $$\operatorname{dim} M=\depth R-1-(\operatorname{grade} M-1),$$
whence the desired equality follows.   \end{proof}

\begin{remark}
One can observe that the equality in the previous theorem can be rewritten as \( \operatorname{qid}_R M = \dim M + \operatorname{grade} M \), using Theorem \ref{inj}.
\end{remark}

By Bass's Conjecture, a local ring is Cohen-Macaulay if it admits a non-zero Cohen-Macaulay module with finite injective dimension. Again, this statement is not true when quasi-injective dimension is considered in place of injective dimension. In \cite[Corollary 3.8]{Gheibi}, Gheibi proved the following: If $R$ is a local ring and there exists a non-zero \( R \)-module \( M \) with maximal Krull dimension and finite quasi-injective dimension, then $R$ is Cohen-Macaulay. As a corollary, we derive a more general sufficient condition for \( R \) to be Cohen-Macaulay. 

\begin{corollary}\label{corol:criteria}
Let $R$ be a local ring. If there exists a non-zero $R$-module $M$ such that $\operatorname{qid}_R M < \infty$ and $\dim R = \dim M + \operatorname{grade} M$, then $R$ is Cohen-Macaulay.
\end{corollary}

\begin{corollary}
    Let $R$ be a local ring, and let $M$ be a non-zero $R$-module. If $M$ is quasi-perfect and $\operatorname{qid}_R M < \infty$, then $M$ is Cohen-Macaulay.  
\end{corollary}
\begin{proof}
    We have $$\depth R-  \depth M=\operatorname{grade} M=\depth R-\dim M,$$
    where the first equality is due to $M$ being quasi-perfect and to the Auslander-Buchs\-baum formula for quasi-projective dimension, while the other is because of Theorem  \ref{teo:forinj}. Thus, $\depth M=\dim M$, and therefore,  $M$ is Cohen-Macaulay.
\end{proof}

\section{Grade inequalities for modules with finite quasi-projective dimension}\label{section6} 

In \cite{Yassemi}, the respective authors established some grade inequalities for modules with finite projective dimension or finite Gorenstein dimension. In this section, using similar arguments to those developed in \cite{Yassemi}, we provide some grade inequalities for mo\-dules of finite quasi-projective dimension.   The inequality $\operatorname{grade} M \leq \operatorname{qpd}_R M$, established in Theorem \ref{theo4.5}, will play a crucial role in this section.


\begin{theorem}\label{teo:grade}
Let $M$, $N$ and $L$ be non-zero $R$-modules such that $\operatorname{qpd}_R N< \infty$ and $\operatorname{qpd}_R L < \infty$. If $\operatorname{Supp} M \subseteq \operatorname{Supp} L$, then
$$\operatorname{grade} L + \operatorname{grade}(M,L) \leq  \operatorname{grade}(M,N)+\operatorname{qpd}_R N.$$
\end{theorem}
\begin{proof}
 Choose $\mathfrak{p} \in \operatorname{Supp} M$ such that $\operatorname{grade}(M,N)=\operatorname{depth}N_{\mathfrak{p}}$. Note that $\operatorname{qpd}_{R_{\mathfrak{p}}} N_{\mathfrak{p}} \leq \operatorname{qpd}_R N <  \infty$, by \cite[Proposition 3.5(1)]{GJT}. Thus, by the Auslander-Buchsbaum formula for quasi-projective dimension, we have: 
    \begin{align*}
  \operatorname{grade}(M,N) &  = \operatorname{depth} R_{\mathfrak{p}}-\operatorname{qpd}_{R_{\mathfrak{p}}} N_{\mathfrak{p}} \\
    & = \operatorname{depth} L_{\mathfrak{p}} + \operatorname{qpd}_{R_{\mathfrak{p}}} L_{\mathfrak{p}} - \operatorname{qpd}_{R_{\mathfrak{p}}} N_{\mathfrak{p}} 
    \\ 
    & \geq \operatorname{grade}(M,L) + \operatorname{grade} L_{\mathfrak{p}}-\operatorname{qpd}_R N \,\,\, \text{ (by Theorem \ref{theo4.5})} \\
    & \geq \operatorname{grade}(M,L) + \operatorname{grade} L-\operatorname{qpd}_R N.
    \end{align*}
\end{proof}

\begin{corollary}
Let $M$ and $L$ be non-zero $R$-modules such that $\operatorname{qpd}_R L < \infty$. If $\operatorname{Supp} M \subseteq \operatorname{Supp} L$, then
\begin{align*}
    \operatorname{grade}(M,L) + \operatorname{grade} L \leq \operatorname{grade} M.
\end{align*}
\end{corollary}
\begin{proof}
It follows directly by Theorem \ref{teo:grade}, setting $N=R$.
\end{proof}

\begin{corollary}
Let $M$ and $N$ be non-zero $R$-modules such that $\operatorname{qpd}_R N < \infty$. One then has the inequality 
\begin{align*}
    \operatorname{grade} M \leq \operatorname{grade}(M,N) + \operatorname{qpd}_R N.
\end{align*}
\end{corollary}
\begin{proof}
It follows directly by Theorem \ref{teo:grade}, setting $L=R$.
\end{proof}

\begin{corollary}
Let 
$M$ and $N$ be non-zero $R$-modules such that $\operatorname{qpd}_R N < \infty$ and $\operatorname{Supp} M \subseteq \operatorname{Supp} N$. Then
\begin{align*}
\operatorname{grade} (M,N) + \operatorname{grade} N \leq \operatorname{grade} M \leq \operatorname{grade}(M,N)+ \operatorname{qpd}_R N.
\end{align*}
In particular, if $N$ is quasi-perfect, then $\operatorname{grade}(M,N)=\operatorname{grade} M - \operatorname{grade}N$.
\end{corollary}
\begin{theorem}
Let $R$ be a local ring, and let $M$ and $N$ be non-zero $R$-modules such that $\operatorname{qpd}_R N < \infty$ and $\operatorname{P}_R (M,N)=0$. Then, for any  $R$-module $L$, we have:
\begin{enumerate}
    \item $\operatorname{grade}(L,\operatorname{Hom}_R(M,N))+ \operatorname{qpd}_R N \geq \operatorname{grade} L$. 
    \item If $\operatorname{Supp} L \subseteq \operatorname{Supp}(\operatorname{Hom}_R (M,N))$, then $\operatorname{grade}L \geq \operatorname{grade}(L,\operatorname{Hom}_R (M,N))+\operatorname{grade}N$.
\end{enumerate}
In particular, if $\operatorname{Supp} L \subseteq \operatorname{Supp}(\operatorname{Hom}_R (M,N))$ and $N$ is quasi-perfect, then the equalities hold.
\end{theorem}
\begin{proof} (1) Choose $\mathfrak{p} \in \operatorname{Supp}L$ such that $$\operatorname{grade}(L,\operatorname{Hom}_R (M,N))=\operatorname{depth}( \operatorname{Hom}_{R_\mathfrak{p}}(M_\mathfrak{p},N_\mathfrak{p})).$$ Since $\operatorname{P}_R(M,N)=0$, then it is easy to see that $\operatorname{P}_{R_{\mathfrak{p}}}(M_{\mathfrak{p}},N_{\mathfrak{p}})=0$. Thus, by \cite[Lemma 4.1]{ArayaYoshino}, we have that $$\operatorname{grade}(L,\operatorname{Hom}_R (M,N))=\operatorname{depth} (\operatorname{Hom}_{R_\mathfrak{p}}(M_\mathfrak{p},N_\mathfrak{p}))=\operatorname{depth} N_\mathfrak{p}.$$
Again, note that $\operatorname{qpd}_{R_{\mathfrak{p}}} N_{\mathfrak{p}} \leq \operatorname{qpd}_R N <  \infty$ by \cite[Proposition 3.5(1)]{GJT}. Therefore, by the Auslander-Buchsbaum formula for quasi-projective dimension, we have:
\begin{align*}
\operatorname{grade}(L,\operatorname{Hom}_R (M,N)) 
& = \operatorname{depth} R_{\mathfrak{p}}- \operatorname{qpd}_{R_\mathfrak{p}} N_{\mathfrak{p}} \\
& \geq \operatorname{grade} L  - \operatorname{qpd}_R N. 
\end{align*}

(2) Choose $\mathfrak{p} \in \operatorname{Supp} L$ such that $\operatorname{grade} L = \operatorname{depth} R_{\mathfrak{p}}$. Again, since $\operatorname{P}_R(M,N)=0$, we have that $\operatorname{depth} (\operatorname{Hom}_{R_\mathfrak{p}}(M_\mathfrak{p},N_\mathfrak{p}))=\operatorname{depth} N_\mathfrak{p}$. Therefore, by the Auslander-Buchsbaum formula for quasi-projective dimension, we have:
\begin{align*}
\operatorname{grade} L  & = \operatorname{depth} N_{\mathfrak{p}} + \operatorname{qpd}_\mathfrak{p} N_\mathfrak{p} \\
& = \operatorname{depth}( \operatorname{Hom}_{R_\mathfrak{p}}(M_\mathfrak{p},N_\mathfrak{p})) + \operatorname{qpd}_{R_{\mathfrak{p}}} N_{\mathfrak{p}} \\
& \geq \operatorname{grade}(L,\operatorname{Hom}_R(M,N)) + \operatorname{grade} N_{\mathfrak{p}} \; \text{ (by Theorem \ref{theo4.5})} \\
& \geq \operatorname{grade}(L,\operatorname{Hom}_R(M,N)) + \operatorname{grade} N.
\end{align*}
\end{proof}
 In the next theorem, which is the main result of this section, we present a refined version of the second inequality of \cite[Theorem 3.1]{ArayaYoshino}, replacing the projective dimension with the quasi-projective dimension. Furthermore, we recover the first inequality by introducing an additional hypothesis.
\begin{theorem}\label{genary}
Let $R$ be a local ring, and let $M$ and $N$ be non-zero $R$-modules such that $\operatorname{qpd}_R N< \infty $ and $\operatorname{q}^R(M,N)=0$. Then, for any $R$-module $L$, we have:
\begin{enumerate}
    \item $\operatorname{grade}(L,M) \leq \operatorname{grade}(L,M \otimes_R N)+ \operatorname{qpd}_R N $.
    \item If $\operatorname{Supp} L \subseteq \operatorname{Supp} N$, then $\operatorname{grade}(L,M \otimes_R N) + \operatorname{grade} N \leq \operatorname{grade}(L,M)$.
\end{enumerate}
In particular, if $\operatorname{Supp} L \subseteq \operatorname{Supp} N$ and $N$ is quasi-perfect, then the equalities hold.
\end{theorem}
\begin{proof}
(1) Choose $\mathfrak{p} \in \operatorname{Supp} L$ such that $\operatorname{grade}(L,M \otimes_R N) =\operatorname{depth}(M \otimes_R N)_{\mathfrak{p}}$. As $\operatorname{q}^R(M,N)=0$, then it is easy to see that $\operatorname{q}^{R_{\mathfrak{p}}}(M_{\mathfrak{p}},N_{\mathfrak{p}})=0$. Therefore, by the Auslander-Buchsbaum formula and the depth formula for quasi-projective dimension (see Theorems \ref{depthformula} and \ref{rem:ABF}) we have that: 
     \begin{align*}
         \operatorname{grade}(L,M \otimes_R N) &  = \operatorname{depth}(M \otimes_R N)_{\mathfrak{p}} \\ 
         & = \operatorname{depth}(M_{\mathfrak{p}} \otimes_{R_{\mathfrak{p}}} N_{\mathfrak{p}})\\
         & = \operatorname{depth} M_{\mathfrak{p}}-\operatorname{qpd}_{R_{\mathfrak{p}}} N_{\mathfrak{p}} \\
  & \geq \operatorname{grade}(L,M)- \operatorname{qpd}_R N.  
     \end{align*}
(2) Choose $\mathfrak{p} \in \operatorname{Supp} L $ such that $\operatorname{grade}(L,M) = \operatorname{depth} M_{\mathfrak{p}}$. Again, we have that as $\operatorname{q}^R(M,N)=0$, then we see that $\operatorname{q}^{R_{\mathfrak{p}}}(M_{\mathfrak{p}},N_{\mathfrak{p}})=0$. Another time, by the Auslander-Buchsbaum formula and the depth formula for quasi-projective dimension, we have that:
\begin{align*}
\operatorname{grade}(L,M) & =  \operatorname{depth} (M_{\mathfrak{p}} \otimes_{R_{\mathfrak{p}}} N_{\mathfrak{p}}) + \operatorname{qpd}_{R_{\mathfrak{p}}} N_{\mathfrak{p}}.\\
& \geq  \operatorname{grade}(L,M \otimes_R N)+ \operatorname{grade} N_{\mathfrak{p}} \; \text{ (by Theorem \ref{theo4.5})} \\
& \geq \operatorname{grade}(L,M\otimes_R N) + \operatorname{grade}N .
\end{align*}
\end{proof}



\begin{fund}
The first author was supported by  S\~ao Paulo Research Foundation (FAPESP) under grant 2019/21181-0. The second author was supported by  S\~ao Paulo Research Foundation (FAPESP) under grant 2022/12114-0. The third author was supported by  S\~ao Paulo Research Foundation (FAPESP) under grants 2022/03372-5 and 2023/15733-5.
\end{fund}
\begin{agra}
Part of this work was developed when the third author was visiting the University of Nebraska-Lincoln. He is grateful for the hospitality. He thanks Roger and Sylvia Wiegand for some discussions. The authors are grateful to the anonymous referee for a careful reading and for valuable comments. 
\end{agra}
\end{document}